\documentclass[12pt]{article}
\usepackage{amsmath,amssymb}
\newtheorem{thm}{Theorem}
\newtheorem{lemma}{Lemma}
\newtheorem{cor}{Corollary}
\newtheorem{prop}{Proposition} 
\newtheorem{rem}{Remark}

\newcommand{\mb}{\mbox}
\newcommand{\x}{\times}
\newcommand{\g}{\mathfrak g}

\newcommand{\ot}{\otimes}
\newcommand{\bs}{\bigskip}
\newcommand{\ga}{\gamma}
\newcommand{\ep}{\epsilon}

\begin{document}
\title{Degenerate Integrability of Spin Calogero-Moser Systems
and the duality with the spin Ruijsenaars systems}

\author{ N.Reshetikhin}
\maketitle

\begin{abstract}
It is shown that spin Calogero-Moser systems are completely
integrable in a sense of degenerate integrability. Their Liouville
tori have dimension less then half of the dimension of the phase
space. It is also shown that rational spin Ruijsenaars systems are
degenerately integrable and
dual to spin Calogero- Moser systems in a sense that action-algle
variables of one are angle-action variables of the other.
\end{abstract}

\section{Introduction}
Let $\mathfrak g$ be a simple complex Lie algebra of rank $r$. Fix
a Borel subalgebra $\mathfrak b$. Let $\mathfrak h$ be the
corresponding Cartan subalgebra of $\mathfrak g$ and $\Delta$ be
the root system of $\mathfrak g$. Choose an enumeration of simple
roots of $\Delta$. Denote by $x_\alpha, \ x_{-\alpha}, \ p_i,
\alpha\in \Delta_+, i=1, \dots, r=rank({\mathfrak g})$ the
coordinates on ${\mathfrak g}^*$ corresponding to the Chevalley
basis $H_i, X_\alpha$ in $\mathfrak g$ ( these coordinates correspond to positive
roots, negative roots and simple roots considered as elements of
the Cartan subalgebra respectively) . The functions
$x_{\alpha_i},x_{-\alpha_i}, p_i$ generate the Poisson algebra of
polynomial functions on ${\mathfrak g}^*$ with the following
determining relations:

\begin{eqnarray*}
\{p_i,p_j\} &=& 0 \\ \{p_i,x_{\alpha_j}\} &=& a_{ij}x_{\alpha_j}
\\ \{p_i,x_{-\alpha_j}\} &=& -a_{ij}x_{-\alpha_j} \\
\{x_{\alpha_i},x_{-\alpha_j}\}&=&\delta_{ij} p_i \\
\{x_{\alpha_i},\dots,\{x_{\alpha_i},x_{\alpha_j}\dots\}&=& 0 ,
i\neq j \\
\{x_{-\alpha_i},\dots,\{x_{-\alpha_i},x_{-\alpha_j}\dots\}&=& 0  ,
i\neq j
\end{eqnarray*}

Let $\mathcal O$ be a regular co-adjoint orbit in ${\mathfrak g}^*$
(i.e. an orbit through a regular semisimple element). It has a natural
symplectic structure. Since the co-adjoint action of $H$ is
Hamiltonian the quotient space ${\mathcal O}/Ad^*_H$ is a Poisson
manifold. It naturally decomposes into the product of Poisson manifolds
${\mathfrak h}^*\times{\mathcal O}'$ where ${\mathcal O}'$ is
symplectic with the
symplectic structure obtained via the Hamiltonian reduction( with
respect to $0\in {\mathfrak h}$) and the Poisson structure on ${\mathfrak
h}^*$ is trivial.

Let $H_{reg}$ be the regular part of the Cartan subgroup in $G$.
Spin Calogero-Moser corresponding to the
orbit $\mathcal O$ is the Hamiltonian system on the symplectic
manifold $T^*(H_{reg})\times {\mathcal O}'$.

The (complex
holomorphic) symplectic structure on this manifold is the product of the symplectic
structure on ${\mathcal O}'$ obtained by Hamiltonian reduction and the
standard symplectic structure $T^*(H_{reg})={\mathfrak h}^*\times H_{reg}$.

Let $\gamma_\alpha$ be a coordinate function on $H$
corresponding to the root $\alpha$ and $p_i$ be the coordinate
function on ${\mathfrak h}^*\subset {\mathcal O}/Ad_H$
corresponding to the simple root $\alpha_i$. We have the following
Poisson brackets between these coordinate functions:
\begin{equation} \label{can-bp}
\{p_i,p_j\}=0, \ \ \{\gamma_\alpha,\gamma_\beta \}=0, \ \{p_i,\gamma_\alpha\}=
(\alpha, \alpha_i) \gamma_\alpha
\end{equation}
where $(.,.)$ is the Killing form on $\mathfrak g$ restricted $\mathfrak h$
and normalized as below.

The Hamiltonian of the spin Calogero-Moser system  is
the following function on $T^*(H_{reg})\times {\mathcal O}'$:

\begin{equation} \label{CM}
H_{CM}=(p,p)/2+\sum_{\alpha \in \Delta_+}\frac{(\alpha,
\alpha)\mu_\alpha
\mu_{-\alpha}}{(\gamma_{\alpha/2}-\gamma_{\alpha/2}^{-1})^2}
\end{equation}

Here $(.,.)$ is the bilinear form on ${\mathfrak h}^*$ induced by
the Killing form $(p,p)= \sum_{ij}p_ip_j (b^{-1})_{ij}$ where
$b_{ij}$ is the symmetrized Cartan matrix, \  $p_i$ and $\gamma_\alpha$
are coordinates on $T^*H$ described above and the product
$\mu_\alpha \mu_{-\alpha}$ is a function on ${\mathcal O}'$ which
is the product of coordinate functions $\mu_\alpha$ and
$\mu_\alpha$ on $\mathcal O$ corresponding to roots $\alpha$ and
$-\alpha$ respectfully.

The normalizer $N(H)\subset G$ of $H$ acts by conjugations on
$\mathcal O$. This induces the natural action of the Weyl group
$W$ on the quotient ${\mathcal O}/Ad^*_H$. The natural Poisson structure on
this space is $W$-invariant. This action of $W$ descends to
the reduced space ${\mathcal O}'$. It is clear that the Hamiltonian of the spin
Calogero-Moser system is
invariant with respect to the diagonal action of $W$ on
$T^*(H_{reg})\x {\mathcal O}'$. Thus, effectively the phase space
of the Calogero-Moser system is
$(T^*(H_{reg})\x {\mathcal O}')/W$. We will call it effective phase space.
This variety is not smooth. Desingualization of this variety is an interseting problem,
which was studied a lot in the case of "usual"( non-spin)  Calogero-Moser systems
in \cite{BW} \cite{EG} (see also references therein).

The complex holomorphic Hamiltonian system described above has two
natural real forms.

\begin{itemize}
\item Compact real form correspond to $p_i$
being real, $|\gamma_\alpha|=1$  and $\mu_\alpha=-\bar{\mu}_{-\alpha}$.

\item Split real form correspond to real $\ p_i, \
\mu_\alpha$ and $\mu_{-\alpha}$ and to $\gamma_{\alpha/2}>0$.
In this case the phase space is
$(T^*{\mathbb R}^r_{>0}\times {\mathcal O}')/W$. It is naturally isomorphic to
$\Delta\times {\mathbb R}^r_{>0}\times{\mathcal O}'$, where $\Delta\subset {\mathbb R}^r$
is the fundamental domain of the Weyl group $W$.
\end{itemize}


Rational spin Calogero-Moser systems (for $sl_n$) were introduced
and studied in \cite{Gibb}. The system (\ref{CM}) and its elliptic
generalization has been studied for $sl_n$ in \cite{Babel-1} where
its complete integrability has been established. It has been
noticed \cite{Babel-2} and \cite{LiXu} that using dynamical $r$- matrices
one can construct Poisson commuting Hamiltonians for spin Calogero-Moser system.

In this note we show that spin Calogero-Moser systems are completely integrable
for all co-adjoint orbits.
More precisely, they are degenerately integrable
\cite{N}: Liouville tori have dimension
$r=rank({\mathfrak g})$, which is less then half of the dimension
of the phase space  (except for special orbits in $sl_n$ case when
the system becomes the usual Calogero-Moser system). Degenerate
integrability occurs in other integrable system related to simple
Lie algebras \cite{Gekh}\cite{Res}.

I would like to thank P. Michor for stimulating discussions
during my stay at ESI and for bringing my attention to the
problem of integrability of spin Calogero-Moser systems.
I am also grateful to M. Gekhtman, B.Kostant and M. Yakimov for
valuable discussions and comments. This work was partly
supported by the NSF grant DMS-0070931 and by the CRDF grant RM-1-2244.
The paper was
completed when the author visited IHES and ENS. The author thanks both
institutions for the hospitality. The main results were announced in
January 2001 at the meeting "Integrability and Quantization" (Toronto).

\section{Degenerate integrability}
The notion of degenerate integrability generalizes the "usual"
Lioville integrability . Degenerately integrable
system of rank $k$ on a $2n$-dimensional symplectic manifold
is a Hamiltonian dynamical system whose trajectories
are parallel to a co-isotropic fibration on the phase space with fibers
of dimension $k<n$.
\subsection{}\label{int}

Let $C({\mathcal M}_{2n})$ be the algebra of functions on a $2n$-dimensional
symplectic manifold ${\mathcal M}_{2n}$. Assume that we have the following
structure on it:
\begin{itemize}
\item $2n-k$ independent functions $J_1,\dots, J_{2n-k}$ functions
which generate Poisson subalgebra $C_J({\mathcal M}_{2n})$.
\item $k$ independent functions $I_1,\dots, I_k$
generating the center of the Poisson
subalgebra $C_J({\mathcal M}_{2n})$.
\end{itemize}

Geometrically, these data means that we have two Poisson projections:
\begin{equation}\label{degint}
{\mathcal M}_{2n}\stackrel{\psi}\longrightarrow B_{2h-k}^J\stackrel{\pi}
\longrightarrow B_k^I
\end{equation}
where $B^J$ and $B^I$ are Poisson varieties (level surfaces of $J_i$ and
$I_i$ respectively), such that connected components of the preimage
$\pi^{-1}(b)\subset B^J$ of a generic
point $b\in B^I$ are symplectic leaves in $B^J$.

Let $H\in C({\mathcal M}_{2n})$ be a function which is in the center
of the Poisson subalgebra $C_J({\mathcal M}_{2n})$ ( in particular, it
Poisson commutes with functions $I_i$). Let  ${\mathcal M}(c_1,\dots,c_{2n-k})=\{x\in
{\mathcal M} | J_i(x)=c_i\}$ be a level surface of functions $J_i$.
We will call it generic relative
to functions $I_1,\dots,I_n$ if the form $dI_1\wedge\dots\wedge
dI_{k}$ does not vanish identically on it. The following is
true \cite{N}:

\begin{thm}\label{deg}
\begin{enumerate}
\item Flow lines of $H$ are parallel to level surfaces of $J_i$.
\item Each connected component of a generic (relative to functions $I_i$)
level surface
has canonical affine structure generated by the flow lines of
$I_1,\dots, I_k$. If the connected component is compact, it is isomorphic
to an $n$-dimensional torus.
\item Flow lines of $H$ are linear in this affine
structure.

\end{enumerate}
\end{thm}

\subsection{}

Let $b\in B^I$ be a generic point and $D$ be  an open neighborhood
of $b$. Choose the
trivialization of ${\pi}$ over $D$:
\[
{\pi}^{-1}(D)\simeq {\pi}^{-1}(b)\x D \ .
\]
Choose a generic  point $c\in B_J$ . Let $c_0\in \pi^{-1}(b)$
be such that with respect to this trivialization $c=(c_0,\{0\})$.
Choose a neighborhood $\tilde{U}\subset {\pi}^{-1}(D)$ of $c$.
The trivialization of $\pi$ gives an isomorphism
$\tilde{U}\simeq D\times U$, where $U\subset{\pi}^{-1}(b)$ is
a neighborhood of $c_0$.  This neighborhood
inherits natural symplectic structure which we will denote $\omega_U$.

Now we can trivialize $\psi$ over $\tilde{U}$. This gives an isomorphism
\[
\psi^{-1}(\tilde{U})\simeq  U\times D\times T .
\]
where $T$ is the fiber of $\psi$ over $c$. If $T$ is compact, it is
isomorphic to a disjoint union of tori.

Functions $I_1,\dots, I_k$ define a local coordinate system on
$D$. Their Hamiltonian flows generate $k$ independent Hamiltonian
vector fields on generic level surfaces of $\psi$. These flows define
affine coordinates $\phi_1, \dots, \phi_k$ on this level surfaces .

\begin{thm} \cite{N72} There exists a trivialization
$f:{\psi}^{-1}(\tilde{U})\simeq  U\times D\times T$  such that
the symplectic form $\omega$ on ${\mathcal M}$ has the form
\[
\omega= f^*\left( \sum^k_{i=1} dI_i\wedge d\varphi_i +
\omega_U)\right)
\]
where $I_i$ are coordinates on $D$, $\phi_i$ are coordinates on $T$ induced
by the Hamiltonian flows of $I_i$ and $\omega_U$ is the symplectic form on the $U$ induced 
by the
Poisson structure on $\pi^{-1}(b)$.
\end{thm}

 The coordinates $I_i,\phi_i$ are called action-angle variables
for degenerate integrable systems. If the Hamiltonian of the system
is a pull-back of a function on $B^I$, its trajectories lie in $T$
and are llines in coordinates $\phi_i$:
\[
\phi_i(t)=\phi_i(0)+\omega_i(H, I)t
\]
where $\omega_i=\frac{\partial H}{\partial I_i}$.

One can replace real smooth manifolds by complex
manifolds (complex algebraic) and Poisson structures by complex
holomorphic (complex algebraic) structures.

\section{The Poisson variety $T^*G//Ad_G$}
\subsection{}Let $G$ be a simple Lie group and
$(.,.)$ be a Killing form on the corresponding Lie algebra
$\mathfrak g$. We will fix its choice throughout this paper.
It fixes a linear
isomorphism $\g\simeq\g^*$.  Let $T^*G$ be the cotangent bundle to
$G$. We can trivialize it by left translations as
$T^*G\simeq{\g}^*\x G$. We will fix this trivialization as well as
the corresponding trivialization of the tangent bundle throughout
the paper.

The adjoint action of the Lie group $G$ on $G$ extends naturally
on $T^*G$. Let $(x,\gamma)\in \g^*\times G$ be a point in the
cotangent bundle,then
\[
h:x\mapsto {\mb{Ad}}^*_h(x) \ , \qquad \gamma\mapsto h\gamma
h^{-1} \ .
\]
This action is Hamiltonian with the moment map
$\mu: T^*G\to {\frak g}^*$
\[
\mu(x,\gamma) =x-{\mb{Ad}}^*_\gamma(x)
\]

Let $T^*G/\!/{\mb{Ad}} \, G$ be the corresponding categorical
quotient ( the spectrum of the ring of $G$-invariant functions on
$T^* G$). Since the action of $G$ is Hamiltonian we have
\begin{prop} The variety $T^*G/\!/{\rm{Ad}} \, G$
has natural Poisson structure.
\end{prop}


Symplectic leaves of $T^*G/\!/{\mb{Ad}} \, G$ can be described
using the moment map.\\ Let ${\mathcal O}_{\mu}$ be the
coadjoint orbit in ${\g}^*$ passing through generic semisimple
point $\mu$.\\ The symplectic leaf in $T^*G/\!/{\mb{Ad}}
\, G$ passing through $[(x,\gamma)]$ is $\mu^{-1}({\mathcal
O}_{\mu(x,\gamma)})//Ad_G$. Its dimension is
\[
\dim(S^{[(x,\gamma)]}) =
\dim(T^*G/\!/{\mb{Ad}}_G)-\dim(G_{\mu(x,\gamma)})= \dim({\mathfrak
g})-\dim(G_{\mu(x,\gamma)})=\dim({\mathcal O})
\]
where $G_a\subset G$ is the stabilizer of $a\in{\g}^*$ with respect to the
Ad${}^*_G$-action. We will denote the symplectic leaf in $T^*G//Ad_G$ corresponding to
the co-adjoint orbit $\mathcal O$ by $S_{\mathcal O}$.

\subsection{}\label{str-sl}Here we will describe the effective phase space
of the spin Callogero-Moser system corresponding to the coadjoint orbit
$\mathcal O$ as an open dense subvariety in $S_{\mathcal O}$. Let
$G_{reg}$ be the subset of regular elements in $G$.

Consider the intersection $\mu^{-1}({\mathcal O})_{reg}={\mathfrak  g}^*\times G_{reg}\cap
\mu^{-1}({\mathcal O})$. 
It  consists of points $(x,\gamma)\in {\mathfrak g}^*\times G_{reg}$
such that
\[
x-Ad_\gamma^* (x)=\mu \in {\mathcal O}.
\]
By the adjoint $G$-action we can always bring a regular $\gamma$ to the form
$\gamma\in H_{reg}\subset G_{reg}$. Assume that $\gamma$ is such.
Then for Chevalley
coordinates of $x$ and $\mu$ we have:
\[
(1-\gamma_\alpha)x_\alpha=\mu_\alpha, \ \alpha\in \Delta
\]
\[
\mu_i=0, \ i=1,\dots, r.
\]
Cartan coordinates of $\mu$ vanish:
\[
\mu_i=0, \ i=1,\dots, r
\]

Thus, on $\mu^{-1}({\mathcal O})_{reg}$ we can choose coordinates
$\mu_{\alpha}, \alpha\in \Delta$,
$\gamma_i$ and $p_i, \  i=1,\dots, r$. It is clear that
$\mu^{-1}({\mathcal O})_{reg}\subset \mu^{-1}({\mathcal O})$ is
invariant with respect to the adjoint $G$-action.

Define the regular part of the symplectic leaf $S_{\mathcal
O}=\mu^{-1}({\mathcal O})//Ad_G$ as $(S_{\mathcal
O})_{reg}=\mu^{-1}({\mathcal
O})_{reg}//Ad_G$. It is clear that $(S_{\mathcal O})_{reg}$ is
open dense in $S_{\mathcal O}$. Taking the corresponding cross-section of
the $Ad_G$-orbit in $T^*G$ we obtain the following statement.
\begin{thm}We have the isomorphism of symplectic varieties
\[
(S_{\mathcal O})_{reg}\simeq (T^*H_{reg}\times {\mathcal O}')/W
\]
where the Weyl group acts diagonally on factors and its action on
${\mathcal O}'$ is described in the introduction.
\end{thm}

Let us describe the coordinate ring $C({\mathcal O}')$. As the
$Ad_G$-module the ring $Pol({\mathfrak g}^*)={\Bbb C}[x_\alpha,
p_i]$ of polynomial functions on ${\mathfrak g}^*$ has the weight
decomposition with $wt(x_\alpha)=\alpha, \ wt(p_i)=0$ and
$wt(ab)=wt(a)+wt(b)$.

\begin{prop}
\begin{itemize}
\item Polynomials of zero weight form the Poisson
subalgebra ${\rm{Pol}} _0({\g}^*)\subset {\rm{Pol}}({\g}^*)$.
\item Functions $\{p_i\}^r_{i=1}$ and
${\rm{Ad}}{}^*_G$-invariant functions generate the center of this
Poisson subalgebra.
\item There is an isomorphism of Poisson algebra ${\rm{Pol}}_0({\g}^*)\simeq {\rm{Pol}}
({\mathfrak
h}^*)\ot_{\mathbb C} A$, where $A$ is the commutative ring
generated by monomials $x_{\beta_1}\dots x_{\beta_n}$ with
$\sum_i\beta_i=0, \ \beta_i\in \Delta$.The Poisson structure on
the ring $Pol({\mathfrak h}^*)$ is trivial and the ring $A$ is isomorphic
to $C({\mathcal O}')$ as a Poisson algebra.
\end{itemize}

\end{prop}

\section{The integrability of spin Calogero-Moser systems}

In this section we construct the pair of Poisson maps as in
(\ref{degint}) in such a way that the Hamiltonians of  spin Calogero-
Moser systems (\ref{CM}) will be pull-backs from functions on the base $B_I$.
According to the section \ref{int} this will prove the integrability of
the system.

\subsection{} Consider the map
\[
\psi': T^*G\simeq{\g}^*\x G\to {\g}^*\x{\g}^*
\]
\[
\psi'(x,\gamma)= (x,-{\mb{Ad}}_\gamma^*(x))
\]
This map is the product of two moment maps, $\mu_{L,R}: T^*G\to {\mathfrak g}^*$
$\mu_L(x,\gamma)=x, \ \mu_R(x,\gamma)=-Ad^*_\gamma(x)$
for left and right action of $G$ on $T^*G$ respectively.

It is cleat that this map is invariant with respect to the natural
$G$-action and it induces the map
\[
\psi: T^*G /\!/{\mb{Ad}}_G \to ({\g}^*\x {\g}^*)/\!/{\mb{Ad}^*}_G
\]
which bring {\rm{Ad}}${}_G(x,\gamma)$ to
{\rm{Ad}}${}^*_G(x,{\rm{Ad}}^*_\gamma(x))$ assuming that $G$ acts
diagonally ${\g}^*\x{\g}^*$ and via co-adjoint action on each of
the factors.

Consider projections
\[
{\tilde\pi}_{1,2}:{\g}^*\x{\g}^*\to {\g}^*
\]
to the first and second factor respectively. They are
Poisson maps and they define Poisson projections
\[
\pi_{1,2}:({\g}^*\x{\g}^*)/\!/{\mb{Ad}}_G \to
{\g}^*/\!/{\mb{Ad}^*}_G \simeq {\mathfrak h}^*/W \ .
\]
where $\mathfrak h$ is a Cartan subalgebra of $\g$. The Poisson
structure on ${\mathfrak g}^*//Ad^*_G$ and trivial. It is also
clear that
\[
\pi_1\circ\psi =\pi_2\circ\psi \ .
\]
Denote the image of $\psi$ in
$({\g}^*\x{\g}^*)/\!/{\mb{Ad}}_G$ by ${\mathcal P}_{\psi}$. Thus we have maps
\[
T^*G/\!/{\mb{Ad}}_G \stackrel{\psi}{\longrightarrow} {\mathcal
P}_{\psi} \stackrel{\pi}{\longrightarrow}
{\g}^*/\!/{\mb{Ad}}_G^*\simeq {\mathfrak h}^*/W
\]
acting as $[(x,\gamma)]\stackrel{\psi}\mapsto [(x,Ad^*_{\gamma}(x))]\stackrel{\pi}
\mapsto [x]$.
Here $[A]$ is the corresponding $G$-orbit passing through $A$.

\subsection{} Let $S_{[(x,\gamma)]}\subset T^*G/\!/{\mb{Ad}}_G$
be the symplectic leaf in $T^*G//Ad_G$ through $[(x,\gamma)]$. As
we have seen in section \ref{str-sl} for the regular part of this symplectic leaf
we have the isomorphism of symplectic varieties $(S_{[(x,\gamma)]})_{reg}\simeq
(T^*H\otimes {\mathcal O}'_{\mu(x,\gamma)})/W$ where ${\mathcal O}'_{\mu(x,\gamma)}$
is the result of the Hamiltonian reduction of the co-adjoint orbit
passing through $\mu(x,\gamma)$ with respect to the co-adjoint action of
$H$ on this orbit.

Restrict the map $\psi$ to the symplectic leaf. We will have
\begin{equation}\label{CM-integr}
S_{[(x,\gamma)]} \stackrel{\psi}{\longrightarrow}
S_{\psi}^{[(x,\gamma)]} \stackrel{\pi}{\longrightarrow}
{\g}^*/\!/{\mb{Ad}}_G^*\simeq {\mathfrak h}^*/W \ .
\end{equation}
Here $S_{\psi}^{[(x,\gamma)]}$ is the image of
$S_{[(x,\gamma)]}$ in $({\g}^*\x{\g}^*)//{\mb{Ad}}_G^*$
with respect to the map $\psi$.

\begin{thm} Let $x\in{\g}^*$ be regular, then,
\[
\dim(\psi^{-1}(\psi([(x,\gamma)])) =r\ .
\]
\end{thm}

{\it Proof.} Assume that
$\psi([(x,\gamma)])=\psi([(x',\gamma')])$, i.e.,
$x={\mb{Ad}}^*_y(x')$, \ ${\mb{Ad}}^*_\gamma(x)={\mb{Ad}}^*_y
{\mb{Ad}}^*_{\gamma'}(x')$ for some $y\in G$. Therefore
\[
x={\mb{Ad}}^*_{\gamma^{-1}y\gamma'y^{-1}}(x) \ ,
\]
which means $\gamma^{-1}y\gamma'y^{-1}=z\in Z_x\subset G$ where
$Z_x$ is the stabilizer of $x\in \g^*$. For $\gamma'$ we have
\[
\gamma'= y^{-1}\gamma z y
\]
which shows that the set of elements $[x',\gamma']$ such that
$\psi[(x',\gamma')]=\psi[(x,\gamma)]$ has the same dimension as
dim$(Z_x)=r$.

Let $[h_x]\in {\mathfrak h}^*/W$ be the image of $x\in{\g}^*$ in
${\g}^*/\!/{\mb{Ad}}_G^*\simeq {\mathfrak h}^*/W$.

\begin{thm} For regular $x\in{\g}^*$ the preimage
$\pi^{-1}([h_x])\subset S_\psi^{[(x,\gamma)]}$ is the \\
symplectic leaf through $\psi([(x,\gamma)])$.
\end{thm}

This theorem follows from the fact that the center of the Poisson algebra $C(T^*G)$ 
coincides
with the pull-back image of the moment map of $Ad^*_G$-invariant
functions on ${\mathfrak g}^*$.

\begin{cor}Symplectic leaves of ${\mathfrak g}^*\times {\mathfrak g}^*$ are
direct products of coadjoint orbits. The diagonal action of $G$ is
Hamiltonian with the moment map $\mu: {\mathfrak g}^*\times
{\mathfrak g}^*\to {\mathfrak g}^*, \ \mu(x,y)=x+y $. Thus, each
symplectic leaf of $({\mathfrak g}^*\times {\mathfrak g}^*)/Ad_G$
is the intersection of $\mu^{-1}$(a coadgoint orbit) with the
product of two coadjoint orbits.
\end{cor}

The image of $\psi$ consists of elements $(x,y)\in {\mathfrak g}^*
\times {\mathfrak g}^*$ satisfying $Ad^*_G(x)=Ad^*_G(y)$. Thus,
symplectic leaves of $\psi(T^*G)$ are intersections of preimages
of coadjoint orbits with respect to the moment map and preimages
of points in ${\mathfrak g^*}/Ad^*_G$ with respect to the map
$\pi$. Thus $\pi^{-1}([h_x])$ is a symplectic leaf in
$S_\psi^{[(x,\gamma)]}$.

\begin{cor} A Hamiltonian system on $(S_{[(x,\gamma)]})_{reg}$ generated by the
pull-back of an {\rm{Ad}}${}_G$-invariant function on ${\mathfrak
g}^*$ is integrable (degenerately). The dimension of its Liouville
tori is $r={\rm{rank}}(G)$.
\end{cor}


\subsection{} Let us derive the Hamiltonian (\ref{CM}) as the pull-back
of a function on ${\mathfrak g}^*//Ad^*_G$ with respect to the maps
in (\ref{CM-integr}). Consider the function
\begin{equation} \label{H}
H(x,\gamma)=\frac{1}{2} (x,x), x\in {\mathfrak g}^*, \gamma\in G
\end{equation}
on $T^*G\simeq{\mathfrak g}^*\times G$ (on the cotangent bundle to $G$ trivialized
by left translations).

Consider its pull-back in the algebra of functions on the symplectic leaf
$S_{\mathcal O}=\mu^{-1}({\mathcal O}//Ad_G$. Consider the restriction of
this function to $(S_{\mathcal O})_{reg}$. The ring of functions on $(S_{\mathcal O})$
is generated by coordinate functions $p_i$, $\gamma_i$ with $\gamma_\alpha\neq 0$
and zero weight monomials in $\mu_\alpha$.
Projecting $(x,\gamma)\in T^*G$ on this intersection we have
\[
x=\sum_ip_iH_i+\sum_{\alpha\in
\Delta}(1-\gamma_\alpha)^{-1}\mu_\alpha X_\alpha
\]
Here $H_i, X_\alpha$ is the Chevalley basis in the Lie algebra $\mathfrak
g$ ( we assume the identification of $\mathfrak g$ with its dual by the
Killing form). Substitute this expression into (\ref{H}) and we obtain the
Hamiltonian of the Calogero-Moser system. This, together with the
previous section proves the integrability of the Caligero-Moser system
corresponding to the co-adjoint orbit $\mathcal O$.

\subsection{}

Consider (here we essentially follow [KKS]) a
co-adjoint orbit through a semisimple element
of minimal nonzero dimension for ${\mathfrak g}=sl_n$.
If $\mu_{ij}$ are matrix coordinates
on $sl_n$ such orbit can be described by coordinates
$\phi_i,\psi_i$, $i=1,\dots ,n$ with

\[
\mu_{ij} =\phi_i\psi_j -\delta_{ij} \frac 1n<\phi,\psi>
\]
where $<\phi,\psi>=\sum^n_{i=1}\phi_i,\psi_i$. All central
functions on such orbits are polynomials in the quadratic Casimir

\[
c=\sum_{ij} \mu_{ij}\mu_{ji} =<\phi,\psi>^2-2
\frac{<\phi,\psi>^2}{n} +\frac 1n <\phi,\psi>^2 =\frac{n-1}{n}
<\phi,\psi>^2 \ .
\]

The constraint $\mu_{ii}=0$ gives $\phi_i\psi_i=
\frac{<\phi,\psi>}{n}$ and therefore the reduced manifold ${\mathcal O}'$ is zero
dimensional with

\[
\mu_{ij}\mu_{ji} =<\phi_i,\psi_i><\phi_j\psi_j>=
\frac{<\phi,\psi>^2}{n^2} =\frac{c}{n(n-1)} \ .
\]

Thus, for such orbits in $sl^*_n$ the Hamiltonian of the spin
Calogero-Moser system becomes the usual Calogero-Moser Hamiltonian

\[
H=<p,p>/2 +\frac{c}{4n(n-1)} \sum_{i<j} \frac{1}{sh(q_i-q_j)^2}
\]
where $p=(p_1,\dots ,p_n)$, $\{p_i,q_i\}=\delta_{ij}$ and $\gamma_{\epsilon_i-
\epsilon_j}=\exp(2q_i-2q_j)$.

The symmetric group $S_n$ ( the Weyl group of $sl_n$) acts on
$(p,q)$ naturally by permutations and the Hamiltonian of
Calogero-Moser system is invariant with respect to this action.

\section{Degenerate integrability of  rational spin Ruijsenaars systems}

A duality relation between spin Calogero-Moser systems and
systems which we will call rational spin Ruijsenaars systems was
observed in \cite{N95} (see also references therein).
We will define rational spin Ruijsenaars
systems below and then will give examples
of such system for $SL_n$. In this section we will show that rational spin
Ruijsenaars systems are degenerately integrable in a very similar
way to spin Calogero-Moser systems, i.e. we will construct a pair of
Poisson projections as in (\ref{degint}) and we will show that
Hamiltonians of rational spin Ruijsenaars systems are pull-backs from
the base of the last projection.

\subsection{}
The algebra of
polynomial functions on $T^*G\simeq{\g}^*\x G$ is the tensor product of
$Pol({\g}^*)\otimes C(G)$ as a commutative algebra. Consider the following
Poisson structure on $T^*G$.

\begin{itemize}
\item The subalgebras $Pol({\g}^*)$ and $C(G)$ are Poisson subalgebras.
\item Poisson brackets between a linear function $X\in{\g}$ on ${\g}^*$
and $f\in C(G)$ is
\[
\{X,f\} = (L_X-R_X)f
\]
where $L_X$ and $R_X$ are left and right invariant vector fields
on $G$ generated by $X$.
\end{itemize}

To distinguish this Poisson structure from the standard
symplectic structure on the cotangent bundle to a manifold
we will write $(T^*G, p)$ for it.

The adjoint action of the group $G$  on $(T^*G,p)$ is Poisson. There
is a natural Poisson map $\tilde\psi':T^*G\to (T^*G,p)$ acting as
$(x,\gamma)\to (x-Ad^*_\gamma(x), \gamma)$.

\begin{rem} We have maps $\psi':T^*G\to {\mathfrak g}^*\times {\mathfrak g}^*$ and
$\tilde{\psi}':T^*G\to {\mathfrak g}^*\times G$. It is interesting
that we have natural isomorphisms $ker(\psi')\simeq
ker(\tilde{\psi)}'\simeq C$ where $C=\{(x,\gamma)\in
T^*G|x=Ad^*_G(x)\}$. This variety plays an important role in
representation theory of real reductive groups \cite{Ko} ( and an
unpublished paper by B.Kostant on a characteristic variety in $T^*G$).
\end{rem}
The map $\psi'$ induces the
Poisson isomorphism
\[
\tilde\psi: T^*G/\!/{\mb{Ad}} \ G \to (T^*G,p)/\!/{\mb{Ad}_G}  \ .
\]
We also have a natural projection
\[
\tilde\pi:(T^*G,p)/\!/{\mb{Ad}} \ G \to G/\!/{\mb{Ad}}  \ G \ .
\]
This projection is Poisson with the trivial Poisson structure on the
base.

Restricting $\tilde\psi$ to a symplectic leaf of
$T^*G/\!/{\mb{Ad}} \ G$ we have
\[
S\stackrel{\tilde\psi}{\mapsto} \tilde\psi(S) \subset
(T^*G,p)/\!/{\mb{Ad}} \ G \stackrel{\tilde\pi}{\to}
G/\!/{\mb{Ad}} \ G \ .
\]

\begin{prop} Over generic point $\dim(\ker\tilde\psi)=r$.\end{prop}

{\bf Proof.} Assume that
$\tilde\psi([(x,\ga)])=\tilde\psi([(x',\ga')])$, i.e., that there
exists $y\in G$ such that
\[
{\mb{Ad}}^*y(\mu(x,\ga)) = \mu(x',\ga') \ , \qquad {\mb{Ad}}(\g) =
\g' \ .
\]
The first equation implies
\[
{\mb{Ad}}^*_y(x)- {\mb{Ad}}^*_y {\mb{Ad}}^*_{\ga}(x) =x' -
{\mb{Ad}}^*_{\ga'}(x)
\]
which, together with the second one, means
\[
x-x'={\mb{Ad}}^*_{\ga'}(x-x') \ .
\]

 From here it follows that the variety of points
$[(x,\ga)]$ such that
$\tilde\psi*([(x,\ga)])=\tilde\psi([x',\ga])$ has dimension $r$.

\bs
It is clear that $\dim(G/\!/{\mb{Ad}}\ G\simeq H/W)=r$ and therefore
we have the pair of Poisson projections
\[
S\stackrel{\tilde\psi}{\longrightarrow} \tilde\psi(S)
  \stackrel{\tilde\pi}{\longrightarrow}
G/\!/{\mb{Ad}} \ G
\]
which satisfy conditions described in section \ref{int}.
In particular \ $\dim(S)=\dim\tilde\psi(S) +\dim(G/\!/{\mb{Ad}}\ G)$.

\begin{prop}
For generic orbit $[h]\in G/\!/{\rm{Ad}}\ G$, \ $\pi^{-1}([h])$ is a
symplectic leaf in $\tilde\psi(S)$.
\end{prop}

{\bf Proof.} The projection of a Hamiltonian flow line on
$(T^*G,p)$ to ${\g}^*$ is parallel to Ad${}^*_G$ action and the
projection of a Hamiltonian flow line to $G$ is parallel to
Ad${}_G$-action. Therefore symplectic leaves in $(T^*G,p)$ are
products ${\mathcal O}\x C$ where ${\mathcal O}\subset{\g}^*$ is a
coadjoint orbit, $C\subset G$ is a conjugacy class. For the same
reason symplectic leaves in $(T^*G,p)/\!/{\mb{Ad}}\ G$ are
${\mathcal O}\x_G C$ and therefore they project to points on $
G/\!/{\mb{Ad}}\ G$.

\begin{cor} Any Hamiltonian flow on a symplectic leaf
$S\subset T^*G/\!/{\rm{Ad}}\ G$ which is generated by the
pull-back of an {\rm{Ad}}${}_G$-invariant function on $G$ is
degenerately integrable.
\end{cor}

Let $\omega_i$ be a fundamental weight of $\mathfrak g$. We will
call functions $H_i=\tilde{\psi}^*\circ\tilde{\pi}^*(\chi_{\omega_i})$
Hamiltonians of rational spin Ruijsenaars systems
because for rank = 1 adjoint orbit in $SL_n$ it coincides with the
rational Ruijsenaars systems [N95](see below). We just proved that all
these Hamiltonians are degenerately integrable.

\subsection{}\label{SL_n}
Consider the special case $G=SL_n$ and the coadjoint orbit of
$rank=1$ . Assume $x\in sl^*_n$ is semisimple and choose a
representative $(h,g)$ of the $SL_n$-orbit through $(x,\ga)$ in
which $h$ diagonal. As in case of spin Calogero-Moser system let
us express $g$ in terms of $x$ and the value of the moment map,
i.e. let us solve the equation
\begin{equation}\label{relation}
(h_i-h_j)g_{ij} = \sum_{k=1}^n\mu_{ik}g_{kj} \ .
\end{equation}
for the non-diagonal elements of $g$.

For $\mu_{ij}$ on this orbit we have:
\[
\mu_{ij}=\phi_i\psi_j-\delta_{ij}\frac{1}{n}<\phi,\psi>
\]

Then the equation (\ref{relation}) can be written as
\[
(h_i-h_j)g_{ij}=\phi_i\sum_k\psi_kg_{kj}-\kappa g_{ij}
\]
where $\kappa=<\phi,\psi>/n$. From this equation we obtain the
identity
\begin{equation}\label{equation-1}
g_{ij}=\frac{1}{h_i-h_j+\kappa}\phi_i\sum_k\psi_kg_{kj}
\end{equation}
This gives the system of equations for $\psi_i\phi_i$
\[
\sum_{i=1}^n\frac{\phi_i\psi_i}{h_i-h_j+\kappa}=1
\]
and the identity
\[
g_{ii}=\frac{\phi_i}{\kappa}\sum_{k=1}^n\psi_kg_{ki}
\]
This, together with the equation (\ref{equation-1}) implies
\[
g_{ij} \ = \ \frac{\phi_i\phi_j^{-1}\kappa g_{jj}}{h_i-h_j+\kappa} \ .
\]
The first two elementary $G$-invariant functions of $g$ are
\begin{eqnarray*}
{\mb{tr}}(g) &=& \sum^n_{i=1}g_{ii} \ , \\ {\mb{tr}}(g^2) &=&
\kappa^2\sum_{ij} g_{ii}g_{jj} \
\frac{1}{(h_i-h_j+\kappa)(h_j-h_i+\kappa)} \ .
\end{eqnarray*}
The reduced Poisson brackets are \
$\{h_i,g_{jj}\}=\delta_{ij}g_{jj}$. The second function is the Hamiltonian
of the rational Ruijsenaars system.

\subsection{}\label{RR} Let us show that in a neighborhood of the
identity in $G$ rational spin Ruijsenaars systems degenerate to rational
Calogero-Moser systems.

First, identify ${\mathfrak g}^*$ with $\mathfrak g$ using the
Killing form. Then assume $x\in \mathfrak g$ is regular and choose a
representative $(h,g)$ of the $G$-orbit in $T^*G$ through $(x,\ga)$ with
$h\in \mathfrak h$. Assume that $g$ is in a neighborhood of $e\in
G$ which is the image of the exponential map and let
$g=\exp(\epsilon)$. Denote by $\epsilon_\alpha, \ep_i$ $\alpha\in
\Delta, i=1,\dots, r$ the Chevalley coordinates of $\ep$.

Coordinates $\ep_\alpha$ can be expressed as a power series in
$\mu$ and $h$. Indeed, we have
\[
h-\exp(ad_\ep)(h)=\mu
\]
It is clear that this equation does not define the Cartan coordinates
$\ep_i$. It is also clear that $\ep_\alpha$ can be expressed as a power series in
monomials $\mu_{\beta_1}\dots\mu_{\beta_1}$ with $\sum_i
\beta_i=\alpha$ such that
\[
\ep_\alpha=\frac{\mu_\alpha}{\alpha(h)}+ O((\frac{\mu}{h})^2)
\]
The Cartan part of $\mu$ can be expressed in terms of $\mu_\alpha$
and $h$.

Thus, any polynomial $Ad_G$ invariant function of $\exp(\ep)$ is a power
series in $\mu$ consisting of the monomials
$\mu_{\beta_1}\dots\mu_{\beta_1}$ with $\sum_i \beta_i=0$.

In particular, for any irreducible finite dimensional
representation $V$ we have
\[
tr_V(g)=dim(V)+ \frac{c_2(V)dim(V)}{dim({\mathfrak g})}
(\frac{1}{2}(\ep^{(c)},\ep^{(c)})+\sum_{\alpha\in \Delta_+}
(\alpha,\alpha)\frac{\mu_\alpha \mu_{-\alpha}}{\alpha(h)^2})
+\dots
\]
where $\ep^{(c)}=\sum_{i=1}^r \ep_i H_i$ is the Cartan part of $\ep$
and $(.,.)_{\mathfrak h}$ is the restriction of the Killing form
to $\mathfrak h$.
The second term in this expression is the Hamiltonian of the
rational spin Calogero-Moser system.

Therefore, in a small neighborhood of the identity in $G$ spin rational
Ruijsenaars system degenerates to the rational Calogero-Moser
system.

\section{Action-angle variables and the duality between rational
spin Ruijsenaars and spin Calogero-Moser systems}

\subsection{Action-angle variables for spin Calogero-Moser systems}

\subsubsection{}The degenerate integrability of spin Calogero-Moser
systems was described earlier by two projections:
\[
S\subset T^*G//Ad_G \stackrel{\psi}{\longrightarrow} \psi(S)\subset ({\mathfrak g}^*\x
{\mathfrak g}^*)//Ad^*_G\stackrel{\pi}{\longrightarrow} {\mathfrak g}^*//Ad^*_G\simeq {\mathfrak
h}^*/W
\]
where $S$ is the phase space of the system, which is a subvariety in a
symplectic leaf of $T^*G//Ad_G$.

Action variables for a spin Calogero-Moser system are
affine coordinates on $Spec(C_G({\mathfrak g}^*))\simeq {\mathfrak
h}^*/W$ (the base of the projection $\pi$). Any linear basis in ${\mathfrak h}^*$
determines such coordinate system. For convenience we will
choose these functions to be dual to simple roots $\alpha_i\in
{\mathfrak h}^*$.

Angle variables, by definition, are affine coordinates on a
fiber over generic point of the projection $\psi$ which are
generated by Hamiltonian flows of pull-backs of action variables.
In these coordinates the Hamiltonian flow generated by the pull-back of a function
on ${\mathfrak h}^*/W$ is linear.

Denote by $F_{(x,\ga)}$ the fiber of $\psi$ which contains the
orbit $Ad_G(x,\ga)$. The following is clear.
\begin{lemma} We have:
\[
F_{(x,\ga)}=Ad_G(x,Z_x\ga)
\]
where $Z_x=\{g\in G| Ad^*_g(x)=x\}$ is the stabilizer of the point
$x$.
\end{lemma}

The Hamiltonian flow generated by the pull-back of the function
$f\in C_G({\mathfrak g})$ on $T^*G$ through the point $(x,\ga)$
has the form
\begin{equation} \label{flow}
(x_t, \ga_t)=(x, \exp(tdf(x)\ga)
\end{equation}
where $df(x)\in \mathfrak g$ is the differential of $f$ at $x\in
{\mathfrak g}^*$. It is clear that $F_{(x,\ga)}$ is invariant with
respect to (\ref{flow}). Points in $F_{(x,\ga)}$ are parameterized
by $Z_x$ as $(x,z\ga)$ where $z\in Z_x$. The Hamiltonian flow
through such point generated by the pull-back $f$ of a function
from $ C_G({\mathfrak g}^*)$ is $z_t=\exp(t df(x)z$.

For regular $x\in {\mathfrak g}^*$ we have the isomorphism of Lie
groups $Z_x\simeq H$. Thus, the same linear coordinate system on
$\mathfrak h$ which gives the action variables, after
exponentiating  gives angle variables on the fiber $F_{(x,\ga)}$
for regular $x$.

\subsubsection{}
Now let us describe how angle variables constructed above can be
expressed in terms of coordinates $p_i, \gamma_\alpha, \mu_\alpha$.
Consider the space of $Ad_G$-orbits in
$T^*G//Ad_G$ through points $(x,\ga)$
with  regular $\gamma$. For the coordinate description of
Hamiltonians of spin Calogero-Moser systems we have chosen the
cross-section   $\{(x,\ga)|x\in {\mathfrak g}, \ \ga\in H_{reg}\}$.

Let us identify $\mathfrak g$ with ${\mathfrak g}^*$ using the
Killing form, then  in coordinates we have used in section 4, for
$(x,\ga)$ we have:
\[
x=\sum^r_{i=1} p_iH_i +\sum_{\alpha\in \Delta}
\frac{\mu_\alpha}{1-\ga_\alpha}X_\alpha, \ \ \ga=(h_{\alpha_1},\dots,
h_{\alpha_r})
\]

If $x$ is semisimple there exists $x_0\in \mathfrak h$ and $u\in
G$ such that
\begin{equation} \label{ss}
x=ux_0u^{-1}
\end{equation}

Assume that the element $u^{-1}\ga u$ has the Gaussian decomposition:
\begin{equation}\label{ss-1}
u^{-1}\ga u =b_+b_-b_0
\end{equation}
where $b_\pm\in N_\pm, \ b_0\in H$.

It is clear that the decomposition (\ref{ss}) is defined modulo the action
\begin{equation} \label{al}
v : (u,x_0)\to (uv, v^{-1} x_0 v )
\end{equation}
of $N(H)$ on pairs $(u,x_0)$. It is also clear that
\begin{itemize}
\item $b_0$ depends only on $Ad_{N(H)}(x,\ga)$,
\item the action (\ref{al}) brings $b_0$ to $v b_0 v^{-1}$,
\item the Hamiltonian flow generated by the pull-back of $f\in C_G({\mathfrak g}^*)\simeq
C_W({\mathfrak h}^*)$ evolve $b_0$ as $b^t_0=b_0 \exp( tdf(x_0))$. Here
$df$ is the differential of $f\in C_W({\mathfrak h}^*)$ regarded as an
element in $\mathfrak h$.
\end{itemize}

Thus, $b_0$ are angle variables for spin Calogero-Moser system and the equations
(\ref{ss}) and (\ref{ss-1}) describe angle variables in terms of initial $(x,\gamma,\mu)$
coordinates.

\bigskip

\bigskip

\subsection{Action-angle variables for Rational Ruijsenaars
systems}

The degenerate integrability of rational spin Ruijsenaars systems
is described by two projections
\[
S\subset T^*G//Ad_G \stackrel{\tilde{\psi}}{\longrightarrow}
 \tilde{\psi}(S)\subset (T^*G,p)//Ad_G\stackrel{\tilde{\pi}}{\longrightarrow}
 G//Ad_G
\]
where $S$ is a symplectic leaf of $T^*G//Ad_G$ which is the phase space of the
system.

The action variables of a rational spin Ruijsenaars system are
affine coordinates on $H/W$. Such coordinates are determined by the
choice of a linear basis in the Lie algebra ${\mathfrak h}=Lie(H)$.

By definition of angle variables, they are affine coordinates on
fibers over generic points of the projection $\tilde{\psi}$ which
are generated by Hamiltonian flows of pull-backs of action
variables.

Denote by $\tilde{F}_{(x,\ga)}$ the fiber of $\tilde{\psi}$ which
contains the orbit $Ad_G(x,\ga)$.

\begin{prop} We have:
\[
\tilde{F}_{(x,\ga)}=Ad_G(x+ C_\ga , \ga)
\]
where $C_\ga=\{x\in {\mathfrak h}^*| Ad^*_\ga(x)=x\}$.
\end{prop}

The Hamiltonian flow generated by the function $f$ which is the pull-back of
a function on $G//Ad_G$ has the form:
\[
(x_t,\ga_t)=(x +t df(\ga), \ga)
\]

For simple $\ga \in G$ we have the isomorphism of vector spaces
$C_\ga\simeq {\mathfrak h}^*$. Thus if we choose action variables
by fixing a linear basis on $\mathfrak h$, the dual basis in
$C_\ga\simeq {\mathfrak h}^*$ defines angle variables for rational
spin Ruijenaars system.

\subsection{The duality between the two systems}

Projections $\psi$ and $\tilde{\psi}$ are dual in a sense that
their fibers meet exactly at one point.

\begin{prop}\label{dual} Let $F_{(x,\ga)}$ and $\tilde{F}_{(x,\ga)}$ be fibers of $\psi$ and
$\tilde{\psi}$ respectively which contain the orbit $Ad_G(x,\ga)$. Then
\[
F_{(x,\ga)}\cap \tilde{F}_{(x,\ga)}=Ad_G(x,\ga)
\]
\end{prop}

Indeed, $F_{(x,\ga)}=\{Ad_G(x,z\ga)|z\in Z_x\}$ and
$\tilde{F}_{(x,\ga)}= \{Ad_G(x+c,\ga)|c\in C_\ga\}$. It is clear
the the intersection of these two sets consists of one point
$Ad_G(x,\ga)$.

Consider projections $p=\pi\circ\psi$ and $\tilde{p}=\tilde{\pi}\circ
\tilde{\psi}$:
\[
\begin{array}{ccc}
S_{\mathcal O}& \stackrel{\tilde{p}}{\longrightarrow}& G//Ad_G
\\ \downarrow p & {} & \\ {\mathfrak g}^*//Ad^*_G&
 {} &
\end{array}
\]

It is clear that for generic $(x,\ga)\in T^*G$ there are birational isomorphisms
\[
p(\tilde{F}_{(x,\ga)})=\{Ad^*_G(x+c)|c\in C_\ga\}\simeq {\mathfrak h}^*/W
\]
\[
\tilde{p}(F_{(x,\ga)})=\{Ad_Gz\ga|z\in Z_x\}\simeq H/W
\]
and therefore angle variables for rational a spin Ruijenaars system
are action variables for the corresponding spin Calogero-Moser system and vice versa.
In this sense the two systems are dual to each other.

\subsection{Self-duality for rational rational spin Calogero-Moser
system}

As it was pointed out in the section \ref{RR} rational spin
Ruijenaars system degenerates to rational Calogero-Moser system
in a neighborhood of the identity in $G$.
The duality between two systems described in the previous section
becomes the self-duality for spin Calogero-Moser system.

Degenerate integrability of rational spin Calogero-Moser systems
follwos from the results of the previous sections. Let us make a
short summary of what how it works. In a small neighborhood of
$e\in G$ the conatngent bundle $T^*G$ can be replaced by
$T^*{\mathfrak g}={\mathfrak g}^*\oplus\mathfrak g$. The the
construction goes as follows:

\begin{itemize}
\item The group $G$ acts naturally on $T^*{\mathfrak g}={\mathfrak g}^*\oplus{\mathfrak g}$
and this action is Hamiltonian.
\item The variety $T^*{\mathfrak g}//G$ is Poisson.
\item Symplectic leaves of $T^*{\mathfrak g}//G$ are $\mu^{-1}({\mathcal O})//G$ where
$\mathcal O$ is a coadjoint orbit and $\mu: T^*{\mathfrak g}\to {\mathfrak g}^*$
is the moment map for the adjoint action of $G$ on $T^*{\mathfrak g}$.
\item Projections $\psi: T^*{\mathfrak g}//G\to ({\mathfrak g}^*\oplus{\mathfrak g}^*)//Ad^*_G$,
$\psi(G(x,a))=G(x,ad^*_a(x))$ and
$\tilde{\psi}: T^*{\mathfrak g}//G\to ({\mathfrak g}\oplus{\mathfrak g}^*)//G$, $\tilde{\psi}(G(x,a))=G(a,ad^*_a(x))$ are Poisson.

\item The images of $\psi$ and $\tilde{\psi}$ are isomorphic as Poisson manifolds.

\item Let $\pi: ({\mathfrak g}^*\oplus{\mathfrak g}^*)//G\to {\mathfrak g}^*//Ad^*_G$
and $\tilde{\pi}: ({\mathfrak g}\oplus{\mathfrak g}^*)//G\to {\mathfrak g}//G$ be
projections to the first summand.

\end{itemize}

Similarly to how it was done above for spin Calogero-Moser systems it is easy
to verify that pais of projections $(\psi, \pi)$ and $(\tilde{\psi}, \tilde{\pi})$
describe degenerate integrability of spin Calogero-Moser systems.

The duality described in the previouse section now becomes the duality between
one copy of rational Calogero-Moser system with another copy of this system.

\section{Conclusion}

As we see now, a spin Calogero-Moser system corresponding to the coajoint
orbit $\mathcal O$ defines an integrable
system on the symplectic manifold $S_{\mathcal O}$ which can be regarded
as a desingularization of $ (T^*H_{reg}\times {\mathcal O})/W$.

If ${\mathcal O}'=\{point\}$ ( this occur for $G=SL_n$ and for a non-nilpotent
coadjoint  orbit of dimension $n-1$) the variety $S_{\mathcal O}$ is the
Hilbert scheme of a point. This variety and its generalizations were
studied in \cite{BW} and \cite{EG}. It seems that less is known about varieties
$S_{\mathcal O}$ for more complicated orbits.

The results of this paper can be easily generalized to "trigonometric" spin
Ruijsenaars system (spin Macdonald-Ruijsenaars) it will be done in the next paper
where we also will show how to quantize such
systems using harmonic analysis on simple Lie groups and
corresponding quantum groups.
Phase spaces of such systems are symplectic leaves of
the moduli space of flat $G$-connections over punctured torus. These systems
are self-dual and the duality map can be regarded as the action of one of the
generators of the modular group of a punctured torus.


\begin{thebibliography}
{99}

\bibitem[GS97]{Gekh}
M.I.~Gekhtman, M.Z.~Shapiro.
\newblock Non-commutative and commutative integrability of generic
Toda flow in simple Lie algberas.
\newblock
{\em Comm. Pure Appl. Math.} 52: 53--84 (1999).

\bibitem[N72]{N}
N.N.~Nekhoroshev.
\newblock Action-angle variables and their generalizations.
\newblock {\em Trans. Moscow Math. Soc.} 26:180-197 (1972).

\bibitem[GH84]{Gibb}
J.~Gibbons, T.~Hermsen.
\newblock A generalization of the Calogero-Moser system.
\newblock {\em Physica} 11D:337-348(1984).

\bibitem[KBBT94]{Babel-1}I.~Krichever, O.~Babelon, E.~Billey and
M.~Talon.
\newblock Spin Generalization of Calogero-Moser system and
the matrix KP equation.
\newblock hep-th/9411160.

\bibitem[ABB96]{Babel-2}
J.Avan, O.Babelon and E.Billey.
\newblock The Gervais-Neveu-Felder Equation and the Quantum
Calogeor-Moser System.
\newblock {\em Comm. in Math. Physics} 178:281-299(1996).

\bibitem[LX00]{LiXu}
L.C.~Li, P.~Xu.
\newblock Spin Calogero-Moser systems associated with simple Lie
algebras
\newblock {\em C.R.Acad. Sci. Paris, Serie I}, 331: n1, 55-61(2000); math.SG/0009180.

\bibitem[KKS78]{KKS}
D. ~Kazhdan, B. ~Kostant and S.~Sternberg.
\newblock Hamiltonian group actions and dynamical systems of
Calogero type.
\newblock {\em Comm. Pure Appl. Mathe} 31:n4, 481-507(1978)

\bibitem[R01]{Res}
N.~Reshetikhin.
\newblock Degenrate integrability of characteristic
Hamiltonian systems.
\newblock math.QA/0103147

\bibitem[N95]{N95}
N.~Nekrasov
\newblock Holomorphic bundles and many-body systems.
\newblock {\it Comm. Math. Phys.}, 180 (1996)587-604.

\bibitem[K]{Ko}
B.~Kostant
\newblock Lie group representations on polynomial rings.
\newblock {\it Amer. J. Math.} 85(1963)327-404.

\bibitem[EG]{EG}
P.~Etingof and V.~Ginzburg
Symplectic reflection algebras, Calogero-Moser space, and deformed 
Harish-Chandra homomorphism.
\newblock preprint, math.AG/0011114.

\bibitem[BW]{BW}
Y.~Berest and G.~Wilson
\newblock
Automorphisms and Ideals of the Weyl Algebra.
\newblock preprint, math.QA/0102190.

\end{thebibliography}
\end{document}